\documentclass[11pt, reqno,oneside]{amsart}

\usepackage{amscd,amsmath}
\usepackage{mathrsfs}
\usepackage{amsfonts}
\usepackage{amssymb}
\usepackage{enumerate}
\usepackage{setspace}
\usepackage{color}
\usepackage{bm}
\usepackage{esint}
\usepackage{float}
\usepackage{tikz}
\usetikzlibrary{positioning,fadings,through}
\usepackage{anysize}
\usepackage[english]{babel}
\usepackage{tabu}
\usepackage[colorlinks=true, citecolor=blue]{hyperref}
\newtheorem{theorem}{Theorem} [section]

\newtheorem{definition}[theorem]{Definition}

\newcommand{\eqn}{\begin{eqnarray}}
\newcommand{\een}{\end{eqnarray}}

\newcommand{\ZZ}{\mathbb{Z}}
\newcommand{\TT}{\mathbb{T}}
\newcommand{\pat}{\partial_t}

\definecolor{luh-dark-blue}{rgb}{0.0, 0.313, 0.608}

\numberwithin{equation}{section}

\usepackage{geometry}
\usepackage{graphicx}

\usepackage{scalerel,stackengine}
\stackMath
\newcommand\reallywidehat[1]{%
\savestack{\tmpbox}{\stretchto{%
  \scaleto{%
    \scalerel*[\widthof{\ensuremath{#1}}]{\kern-.6pt\bigwedge\kern-.6pt}%
    {\rule[-\textheight/2]{1ex}{\textheight}}
  }{\textheight}%
}{0.5ex}}%
\stackon[1pt]{#1}{\tmpbox}%
}
\usepackage{soul}
\usepackage{mathtools}%

\def\ds{\displaystyle}

\newcommand{\m}[1]{{\mathcal{#1}}}

\newcommand{\w}[1]{{\reallywidehat{#1}}}

\def\q{\quad}
\def\qq{\qquad}

\def\vp{\varphi}
\def\pa{\partial}

\def\N{\nabla}
\def\D{\Delta}

\newcommand{\pare}[1]{\left( #1 \right)}
\newcommand{\norm}[1]{\left\| #1 \right\|}
\newcommand{\av}[1]{\left| #1 \right|}

\newcommand{\set}[1]{\left\{ #1 \right\}}

\renewcommand{\t}[1]{\text{#1}}
\newcommand{\dive}[1]{\t{div}\left( #1 \right)}

 \usepackage[new]{old-arrows}

\allowdisplaybreaks

\begin{document}

\setlength{\abovedisplayskip}{5pt}
\setlength{\belowdisplayskip}{5pt}

\setlength{\jot}{7pt}

\title{Global existence for certain fourth order evolution equations}

\author{Rafael Granero-Belinch\'on}
\address{Departamento de Matem\'aticas, Estad\'istica y Computaci\'on, Universidad de Cantabria, Spain.}
\email{rafael.granero@unican.es}

\author{Martina Magliocca}
\address{Departamento de An\'alisis Matem\'atico, Universidad de Sevilla, Spain.}
\email{mmagliocca@us.es}

\date{\today}


\begin{abstract}
In this paper we establish three global in time results for two fourth order nonlinear parabolic equations. The first of such equations involves the Hessian and appears in epitaxial growth. For such equation we give conditions ensuring the global existence of solution. For certain regime of the parameters, our size condition involves the norm in a critical space with respect to the scaling of the equation and improves previous existing results in the literature for this equation. The second of the equations under study is a thin film equation with a porous medium nonlinearity. For this equation we establish conditions leading to the global existence of solution.
\end{abstract}



\maketitle

\section{Introduction and main results}
High order partial differential equations are a very interesting topic due to their many applications such as beam dynamics, thin films \cite{1-54,1-58,1-59,1-65}, crystal dynamics \cite{GM,M,LS} and many others. From a mathematical viewpoint, its study is more challenging that standard second order PDEs for instance due to its lack of maximum principles and some other features.

Our main goal is proving existence and decay results in Wiener spaces to the following problems:
\begin{equation}\label{eq:pb-hessian}
\begin{cases}
\begin{array}{ll}
\ds\pat u
=
K_0\D u+2K_1\det D^2u-K_2\D^2u-\frac{K_3}{2}\D\pare{\D u}^2
 & \text{in}\q  (0,T)\times\TT^2, \\
\ds u(0,x)=u_0(x) &\t{in}\q \TT^2,
\end{array}
\end{cases}
\end{equation}
with
\[
K_0\ge0, \q K_1\ge0,\q K_2>0,\q K_3\ge0,
\]
and
\begin{equation}\label{eq:pb-porous-thin}
\begin{cases}
\begin{array}{ll}
\ds\pat u
=
-\dive{u\N\D u}-\chi\D u^p
 & \text{in}\q  (0,T)\times\TT^2, \\
\ds u(0,x)=u_0(x) &\t{in}\q \TT^2,
\end{array}
\end{cases}
\end{equation}
with
\[
\chi>0.
\]

Equation in \eqref{eq:pb-hessian} models epitaxial growth, and its geometrical derivation can be found in \cite{2,EGP}.
Roughly speaking, this growth process consists in the superposition of layers due to deposition of new material, all under high vacuum conditions. As pointed out in \cite{EGP}, this phenomenon has several applications such as crystal growth. The authors of this work considered the case $K_0=K_3=0$, and they proved the existence of solutions to
\[
\ds\pat u
=
\det D^2u-\D^2u+f(t,x)\qq
 \text{in}\q  (0,T)\times\Omega,\qq \Omega\subsetneq\mathbb{R}^2,
\]
which are global in time under smallness assumptions on the  data, or local in time with arbitrary data. They assumed $H^2$ initial data, and they studied both the homogeneous case $f\equiv0$ and the case $f\in L^2(0,T;L^2(\Omega))$ with Dirichlet or Navier boundary conditions. When $f\equiv0$ and $\norm{u_0}_{H^2}$ is large enough, they also proved that the $W^{1,4}$ norm of the solution blows up in finite time. Later, Escudero \cite{2} improved the  blow up result previously obtained in \cite{EGP}.

\medskip

Another interesting model of which is linked to \eqref{eq:pb-porous-thin} is contained in \cite{3}. Here,  the authors  studied a problem modeling the effect of odd viscosity on the instability of liquid film along a wavy inclined bottom with linear temperature variation, and they found that the free boundary evolution equation verifies the following asymptotic equation:
\[
\pat u =-A(u)\pa_xu-\alpha\pa_x\pare{B(u)\pa_xu+C(u)\pa_x^3u},
\]
which, for appropriate choice of $A,B$ and $C$, is equivalent to \eqref{eq:pb-porous-thin}.

There exists a huge literature concerning equations as in \eqref{eq:pb-porous-thin}, and  exhaustive resume is given in \cite{1}, together with sharp conditions providing the existence of global solutions. Equation
\begin{equation}\label{eq:1}
\ds\pat u
=
-\dive{u^n\N\D u} -\chi\D u^p,\qq u>0,\q \chi\in\mathbb{R},
\end{equation}
describes the evolution of the thin-film liquid height $u$ spreading on a solid surface. The fourth order term takes into account the surface tension, while the porous medium one is related to the Van der Waals forces. We refer to \cite[Section 1.1]{1-7} for further physical explanations.

The case $\chi=0$, $n=1$, has been deeply studied in \cite{1-54}. The author also provides the derivation of the one-dimensional model, and a detailed description of the physical experiment motivating the interest of the equation itself. The cases with $\chi=0$ and $n\in(0,2)$ and $n\in[2,3)$ can be found in \cite{B98,DP01}  and \cite{1-41} respectively.

The case $\chi=-1$ in \eqref{eq:1}
has been dealt with in  \cite{1-28,1-7}, and the relation among the positive parameters $p$ and $n$ has been investigated.

The blow up of solutions has been proved in the one-dimensional case and for $\chi=1$ in \cite{1-8} when $p\ge n+3$, and this result has been refined in \cite{1-65} in the critical case $p=n+3$.\\
Always in the case $\chi=1$, the  existence of self-similar solutions to \eqref{eq:1}, as well as blow up results, are contained in \cite{1-59} with $p=n+3$,  \cite{1-58}  with $0<n<3$, $n\le p$, and \cite{1-38,1-37} for the case of the first critical exponent $p=n+1+2/N$, $0<n<3/2$, $N\ge1$.

Finally, equations as \eqref{eq:pb-porous-thin} are related also to approximations of nonlocal aggregation-diffusion models  \cite{1-5,1-34}, and tumor growth \cite{1-25,1-39}.

\medskip

In this paper we are going to establish the global existence of weak solutions and decay assuming Wiener initial data. The technique we are going to apply is contained in \cite{GM} (see also \cite{M} and \cite{BG1,BG2,LS,A}). We present the notions of definitions of weak solutions we consider below.

\begin{definition}[Weak solution to \eqref{eq:pb-hessian}]\label{def:ex-hessian}
We say that a function $u$ is a weak solution of \eqref{eq:pb-hessian} if
\[
u\in L^\infty( (0,T)\times\TT^2)\cap L^2(0,T;W^{2,4}(\TT^2))
\]
and verifies the following weak formulation:
\begin{align*}
&\int_{\TT^2} u_0\vp(0)\,dx\\
&\q+\iint_{\TT^2\times(0,T)} u\pat\vp+\pare{K_0 u+\frac{K_3}{2}(\D u)^2}\D\vp+2K_1 \vp \det D^2u
-K_2u\D^2\vp  \,dx\,dt=0
\end{align*}
for every
\[
\vp\in W^{1,1}(0,T;L^1(\TT^2))\cap L^1(0,T;W^{4,1}(\TT^2))\cap L^2(0,T;H^2(\TT^2)).
\]
\end{definition}

Taking advantage of the fact that the equation \eqref{eq:pb-porous-thin} conserves the mean, we define the new variable
\[
v(t,x)=u(t,x)-\frac{1}{4\pi^2}\int u_0(x)\,dx.
\]
Without loss of generality from this point onwards we assume that
$$
\frac{1}{4\pi^2}\int u_0(x)=1.
$$
Hence, \eqref{eq:pb-porous-thin} becomes
\begin{equation}\label{eq:pb-porous-thin-v}
\begin{cases}
\begin{array}{ll}
\ds\pat v
=-\D^2v
-\dive{v\N\D v}-\chi\D (1+v)^p
 & \text{in}\q  (0,T)\times\TT^2, \\
\ds v(0,x)=v_0(x)=u_0(x)-1 &\t{in}\q \TT^2.
\end{array}
\end{cases}
\end{equation}

\begin{definition}[Weak solution to \eqref{eq:pb-porous-thin-v}]\label{def:ex-porous-thin}
We say that a function $v$ is a weak solution of \eqref{eq:pb-porous-thin-v} if
\[
v\in L^\infty( (0,T)\times\TT^2)\cap L^2(0,T;H^2(\TT^2)),
\]
and verifies the following weak formulation:
\[
-\int_{\TT^2} v_0\vp(0)\,dx+\iint_{\TT^2\times(0,T)} -v\pat\vp-v\nabla\D v\cdot \nabla \vp+\chi (v+1)^{p-1}\D\vp  \,dx\,dt=0
\]
for every
\[
\vp\in W^{1,1}(0,T;L^1(\TT^2))\cap L^{q'}(0,T;H^2(\TT^2))\q \t{with}\q1\le q<2.
\]
\end{definition}

The $k$-th Fourier coefficients of a $2\pi$-periodic function on $\mathbb{T}^{d}$ are
\[
\widehat{u}(k)=\frac{1}{(2\pi)^{d}}\int_{\mathbb{T}^{d}}u(x)e^{-ik\cdot x}dx,
\]
and the Fourier series expansion of $u$ is given by
\[
u(x)=\sum_{k\in \mathbb{Z}^{d}}\widehat{u}(k) e^{ik\cdot x}.
\]
Using this, we define the Wiener spaces for $s\ge 0$
$$
A^s=\left\{u\in L^{1}(\mathbb{T}^{d}): \|u\|_{\dot{A}^s}=\sum_{k\in \mathbb{Z}^{d}}|k|^{s}\left|\widehat{u}(k)\right|<\infty\right\}.
$$
We note that $A^0$ is a Banach algebra and furthermore,
$$
A^s\subset C^s\subset H^{s}.
$$

Theorems \ref{teo:ex-hessianA2} and \ref{teo:ex-hessianA0} contains two existence results for problem \eqref{eq:pb-hessian}, and the main difference concerns the regularity of the initial data.

\begin{theorem}[Existence result to \eqref{eq:pb-hessian} for $K_3>0$]\label{teo:ex-hessianA2}
Let $K_3>0$ in \eqref{eq:pb-hessian} and consider $u_0\in A^2$ be a zero mean initial data such that
\begin{equation}\label{eq:smallness-hessian}
K_2-2(K_3+K_1)\norm{u_0}_{A^2}>0\qq \t{if }K_0=0,
\end{equation}
and
\[
\min\set{K_2-2K_3\norm{u_0}_{A^2},K_0-2K_1\norm{u_0}_{A^2}}>0\qq \t{if }K_0>0.
\]
Then there exists at least one global weak solution of \eqref{eq:pb-hessian} in the sense of Definition \ref{def:ex-hessian}
$$
u\in L^\infty(0,T;A^2)\cap L^1(0,T;A^6)\quad \forall T.
$$
Furthermore, the solution satisfies
$$
\norm{u(t)}_{A^2}\le e^{-(K_2-2(K_3+K_1)\norm{u_0}_{A^2})t}\norm{u_0}_{A^2}\qq\t{if }K_0=0,
$$
and
$$
\norm{u(t)}_{A^2}\le e^{-\min\set{K_2-2K_3\norm{u_0}_{A^2},K_0-2K_1\norm{u_0}_{A^2}}t}\norm{u_0}_{A^2}\qq\t{if }K_0>0.
$$
\end{theorem}

Similarly, in the case $K_3=0$, we have that:
\begin{theorem}[Existence result to \eqref{eq:pb-hessian} for $K_3=0$]\label{teo:ex-hessianA0}
Let $K_3=0$ in \eqref{eq:pb-hessian} and consider $u_0\in A^0$ be a zero mean initial data such that
\begin{equation*}
K_2-2K_1\norm{u_0}_{A^0}>0.
\end{equation*}
Then there exists at least one global weak solution of \eqref{eq:pb-hessian} in the sense of Definition \ref{def:ex-hessian}
$$
u\in L^\infty(0,T;A^0)\cap L^1(0,T;A^4)\quad \forall T.
$$
Furthermore, the solution satisfies
$$
\norm{u(t)}_{A^0}\le e^{-(K_2-2K_1\norm{u_0}_{A^0})t}\norm{u_0}_{A^0}.
$$
\end{theorem}

This particular result improves the previous global in time result contained in \cite{EGP} due to the fact that our size condition is given in $A^0$ instead of $H^2$. In fact, the space $A^0$ is a critical space with respect to the scaling of the equation
$$
u_\lambda(x,t)=u(\lambda x,\lambda^4 t).
$$

\medskip

We now present our existence result concerning problem \eqref{eq:pb-porous-thin-v}.

\begin{theorem}[Existence result to \eqref{eq:pb-porous-thin-v}]\label{teo:ex-porous}
Let $\chi<1$ and $0\leq u_0\in A^0$ be a an initial data satisfying
$$
\frac{1}{4\pi^2}\int u_0dx=1
$$
together with the smallness condition
\begin{equation*}
1-\chi-2\norm{v_0}_{A^0}
-\frac{c\chi p!}{2}\pare{\norm{v_0}_{A^0}+2
\sum_{q=1}^{p-1}\norm{v_0}_{A^0}^q
}>0.
\end{equation*}
Then there exists at least one global weak solution of \eqref{eq:pb-porous-thin-v} in the sense of Definition \ref{def:ex-porous-thin}
$$
v\in L^\infty(0,T;A^0)\cap L^1(0,T;A^4)\quad \forall T.
$$
Furthermore, the solution satisfies
\begin{equation*}
\norm{v(t)}_{A^0}\le \exp\pare{-\pare{
1-\chi-2\norm{v_0}_{A^0}
-c\chi p!\pare{\norm{v_0}_{A^0}+2
\sum_{q=1}^{p-1}\norm{v_0}_{A^0}^q
}}
t}
\norm{v_0}_{A^0}.
\end{equation*}
\end{theorem}

In the following, we write
\[
f,_j=\pa_{x_j}f
\]
for the space derivative in the $j-$th direction.

\section{Proof of Theorems \ref{teo:ex-hessianA2} and \ref{teo:ex-hessianA0}}
\textbf{Approximate problem:} We explicitly compute each term of \eqref{eq:pb-hessian}:
\begin{align*}
\D u&=u,_{ii},\\
2\det D^2u&=u,_{ii}u,_{jj}-u,_{ij}u,_{ij},\\
\D^2 u&=u,_{iijj},\\
\D\pare{\D u}^2&=(u,_{jj}u,_{kk}),_{ii}=u,_{iijj}u,_{kk}+u,_{jji}u,_{ikk}.
\end{align*}
Hence, \eqref{eq:pb-hessian} is equivalent to
\begin{equation*}
\begin{cases}
\begin{array}{ll}
\ds\pat u
=
K_0u,_{ii}+K_1\pare{u,_{ii}u,_{jj}-u,_{ij}u,_{ij}}-K_2u,_{iijj}-K_3u,_{ii}u,_{jjkk}-K_3u,_{kki}u,_{jji}
 & \text{in}\q  (0,T)\times\TT^2, \\
\ds u(0,x)=u_0(x) &\t{in}\q \TT^2.
\end{array}
\end{cases}
\end{equation*}

We consider the following approximating problem
\begin{equation}\label{eq:pb-hessian-einstein-appr}
\begin{cases}
\begin{array}{lll}
&\ds\pat u^{(n)}
=
K_0u^{(n)},_{ii}+K_1\pare{u^{(n)},_{ii}u^{(n)},_{jj}-u^{(n)},_{ij}u^{(n)},_{ij}}&\\
&\qq\qq-K_2u^{(n)},_{iijj}-K_3(u^{(n)},_{ii}u^{(n)},_{jjkk}+u^{(n)},_{kki}u^{(n)},_{jji})& \text{in}\q  (0,T)\times\TT^2,\\
&\ds u^{(n)}(0,x)=P_n u_0(x) &\t{in}\q \TT^2,
\end{array}
\end{cases}
\end{equation}
where $P_n$ is the projector on the Fourier modes satisfying
$$
|k|\leq n.
$$
The local existence of smooth solutions for such approximate problem can be obtained by using standard methods.

\textbf{A priori estimates in Wiener spaces:} Let us omit the superscript $(n)$ in the following computations. We rewrite each term of the r.h.s. of \eqref{eq:pb-hessian-einstein-appr} in Fourier:
\begin{align*}
\w{u,_{ii}u,_{jj}-u,_{ij}u,_{ij}}(t,k)&=\sum_{m\in\ZZ^2} \pare{\av{m}^2\av{k-m}^2-m_im_j(k_i-m_i)(k_j-m_j)}\times\\
&\qq\times\w{u}(t,m)\w{u}(t,k-m),\\
\w{u,_{iijj}}(t,k)&=\av{k}^4\w{u}(t,k),\\
\w{u,_{ii}}(t,k)&=\av{k}^2\w{u}(t,k),\\
\w{u,_{ii}u,_{jjll}}(t,k)&=-\sum_{m\in\ZZ^2}|m|^4|k-m|^2\w{u}(t,m)\w{u}(t,k-m)\\
\w{u,_{ill}u,_{ijj}}(t,k)&=-\sum_{m\in\ZZ^2}m_i|m|^2(k_i-m_i)|k-m|^2\w{u}(t,m)\w{u}(t,k-m).
\end{align*}
Note that the contribution of the term
$$
K_0\Delta u
$$
is always negative in Wiener spaces.

Assume that $K_0=0$. Then, the Fourier coefficient of  \eqref{eq:pb-hessian-einstein-appr} is given by
\begin{equation*}
\begin{aligned}
\pat \w{u}(k,t)&=K_1\sum_{m\in\ZZ^2} \pare{\av{m}^2\av{k-m}^2-m_im_j(k_i-m_i)(k_j-m_j)}\w{u}(t,m)\w{u}(t,k-m)\\
&\q
-K_2|k|^4\w{u}(t,k)+K_3\sum_{m\in\ZZ^2}(|m|^4|k-m|^2+m_i|m|^2(k_i-m_i)|k-m|^2)\w{u}(t,k-m)\w{u}(t,m).
\end{aligned}
\end{equation*}
We use the fact that
\begin{equation}\label{eq:der-t}
\pat|\w{u}(t,k)|=\frac{Re\left( \overline{\w{u}}(t,k)\pat \w{u}(t,k) \right)}{|\w{u}(t,k)|},
\end{equation}
and also
\[
\sum_{k\in\ZZ^2}	|k|^2\pat\av{ \w{u}(k,t)}=\frac{d}{dt}\norm{u(t)}_{A^2},\qq \sum_{k\in\ZZ^2}|k|^6\av{\w{u}(t,k)}= \norm{u(t)}_{A^6} ,
\]
\begin{align*}
&\sum_{k\in\ZZ^2}|k|^2\av{\sum_{m\in\ZZ^2} \pare{\av{m}^2\av{k-m}^2-m_im_j(k_i-m_i)(k_j-m_j)}\w{u}(t,m)\w{u}(t,k-m)}\\&\le2 \norm{u(t)}_{A^2}\norm{u(t)}_{A^4},
\end{align*}
\begin{align*}
\sum_{k\in\ZZ^2}|k|^2\av{\sum_{m\in\ZZ^2}|m|^4|k-m|^2\w{u}(t,m)\w{u}(t,k-m)}&\le \norm{u(t)}_{A^2}\norm{u(t)}_{A^6},
\end{align*}
\begin{align*}
\sum_{k\in\ZZ^2}|k|^2\av{\sum_{m\in\ZZ^2}m_i|m|^2(k_i-m_i)|k-m|^2\w{u}(t,m)\w{u}(t,k-m)}&\le \norm{u(t)}_{A^3}\norm{u(t)}_{A^5}\leq \norm{u(t)}_{A^2}\norm{u(t)}_{A^6},
\end{align*}
to estimate the  $A^2$ seminorm of $\pat u$ as
\begin{equation}\label{eq:ineq-A2}
\frac{d}{dt}\norm{u(t)}_{A^2}\le -\pare{K_2-2K_3\norm{u(t)}_{A^2}}\norm{u(t)}_{A^6}+2K_1 \norm{u(t)}_{A^2}\norm{u(t)}_{A^4}.
\end{equation}
We now estimate the last term in the r.h.s. of \eqref{eq:ineq-A2} as below
\[
2K_1 \norm{u(t)}_{A^2}\norm{u(t)}_{A^4}\le 2K_1 \norm{u(t)}_{A^2}\norm{u(t)}_{A^6},
\]
so that
\[
\frac{d}{dt}\norm{u(t)}_{A^2}\le -\pare{K_2-(2K_3+2K_1)\norm{u(t)}_{A^2}}\norm{u(t)}_{A^6}.
\]

Thus, if $u_0\in A^2(\TT^2)$ is such that
$$
K_2-2(K_3+K_1)\norm{u_0}_{A^2}>0
$$
and using a contradiction argument in time, we obtain  that $u$ is uniformly bounded in
\begin{equation}\label{eq:unif-bound-hessian}
 W^{1,1}(0,T;A^2(\TT^2))\cap L^1(0,T;A^6(\TT^2)),
\end{equation}
and, furthermore, it decays
\begin{equation*}
\norm{{u}(t)}_{A^2}\le e^{-(K_2-2(K_3+K_1)\norm{u_0}_{A^2})t}\norm{u_0}_{A^2}.
\end{equation*}

If $K_0>0$, we can improve the smallness condition \eqref{eq:smallness-hessian} in the following way. \\
The term $-K_0\av{k}^2\w{u}(t,k)$ appears in the right hand side of the previous equation. Then, reasoning as in the case $K_0=0$, the inequality in \eqref{eq:ineq-A2} takes the following form
\[
\frac{d}{dt}\norm{u(t)}_{A^2}\le -\pare{K_2-2K_3\norm{u(t)}_{A^2}}\norm{u(t)}_{A^6}-\pare{K_0-2K_1 \norm{u(t)}_{A^2}}\norm{u(t)}_{A^4}.
\]
We can thus avoid to estimate $\norm{u(t)}_{A^4}$ with $\norm{u(t)}_{A^6}$ requiring $u_0\in A^2$ such that
\[
\min\set{K_2-2K_3\norm{u_0}_{A^2},K_0-2K_1\norm{u_0}_{A^2}}>0.
\]

\medskip

Moreover, in the case where $K_0=K_3=0$, we can improve the previous result and find that
\begin{equation*}
\frac{d}{dt}\norm{u(t)}_{A^0}\le +2K_1\norm{u(t)}_{A^2}^2-K_2 \norm{u(t)}_{A^4}.
\end{equation*}
Using interpolation, we obtain that
$$
\frac{d}{dt}\norm{u(t)}_{A^0}\le (-K_2+2K_1\norm{u(t)}_{A^0}) \norm{u(t)}_{A^4},
$$
from where we can conclude the desired estimates as before.\\
An analogous estimate holds in the case $K_0>0$.

\textbf{Compactness results:} Up to subsequences, we have that
\begin{align}
u^{(n)}&\to u && \t{a.e. }(0,T)\times \TT^2,\label{eq:ae-hessian}\\
u^{(n)}&\overset{*}{\rightharpoonup} u && \t{in } L^{\infty}(0,T;W^{2,\infty}( \TT^2)),\label{eq:weak*1-hessian}\\
u^{(n)}&\overset{*}{\rightharpoonup} u && \t{in } \m{M}(0,T;W^{6,\infty}( \TT^2)),\label{eq:weak*2-hessian}\\
u^{(n)}&\overset{*}{\rightharpoonup} u && \t{in } L^{\frac{2p}{3p-6}}(0,T;W^{p,\infty}( \TT^2))\q\t{with}\q 2< p< 6,\label{eq:weak-hessian}\\
%
%
u^{(n)}&\to u && \t{in } L^{2}(0,T;H^{2}( \TT^2)).
\label{eq:strong-hessian}
\end{align}
Indeed, the weak-$*$ convergences \eqref{eq:weak*1-hessian} and \eqref{eq:weak*2-hessian} follow from the uniform bound in \eqref{eq:unif-bound-hessian} and the Banach-Alaoglu Theorem. We use the interpolation inequality to say that,
\[
\norm{u^{(n)}(t)-u(t)}_{A^p}\le \norm{u^{(n)}(t)-u(t)}_{A^2}^{\frac{6-p}{2p}}\norm{u^{(n)}(t)-u(t)}_{A^6}^{\frac{3p-6}{2p}},
\]
and we integrate in time
\[
\int_0^T\norm{u^{(n)}(t)-u(t)}_{A^p}^{\frac{2p}{3p-6}}\,dt\le \norm{u^{(n)}-u}_{L^\infty
(A^2)}^{\frac{6-p}{3p-6}}\int_0^T\norm{u^{(n)}(t)-u(t)}_{A^6}\,dt.
\]
Then, \eqref{eq:weak-hessian} follows recalling \eqref{eq:unif-bound-hessian} and observing that  $ \norm{f}_{W^{p,\infty}}\le \norm{f}_{A^p}$. Similarly, the a.e. convergence in \eqref{eq:ae-hessian} follows from the previous ones.

We now focus on the strong convergence in \eqref{eq:strong-hessian}. Interpolation inequality implies that
\[
\norm{u^{(n)}-u}_{L^2(H^2)}\le \norm{u^{(n)}-u}_{L^2(L^2)}^{\frac{1}{3}}\norm{u^{(n)}-u}_{L^2(H^3)}^{\frac{2}{3}}\le c\norm{u^{(n)}-u}_{L^2(L^2)}^{\frac{1}{3}}
\]
by \eqref{eq:weak-hessian} with $p=3$. Then, we have just to prove the strong convergence
\begin{align}\label{eq:strong-L2-hessian}
u^{(n)}&\to u && \t{in } L^{2}(0,T;L^2( \TT^2))
\end{align}
to deduce \eqref{eq:strong-hessian}.
We want to apply classical compactness results (see, for instance, \cite[Corollary 4]{simon}) in order to get \eqref{eq:strong-L2-hessian}. Then, we need some spaces $X$ and $Y$ such that
\begin{align*}
&u^{(n)}\t{ uniformly bounded }L^2(0,T;X),\\
&\pat u^{(n)} \t{ uniformly bounded }L^1(0,T;Y),
\end{align*}
verifying
\[
X\overset{compact}{\varlonghookrightarrow} L^2(\TT^2)\varlonghookrightarrow Y.
\]
We set $X=H^2(\TT^2)$, so the boundedness in $L^2(0,T;H^2(\TT^2))$ follows from the one in $L^2(0,T;H^3(\TT^2))$ and the finiteness of the domain. We choose $Y=H^{-2}(\TT^2)$ and we prove the uniform boundedness of the time derivative using that
\[
\|\pat u^{(n)}(t)\|_{H^{-2}}=\sup_{
\footnotesize
\begin{array}{c}
\vp\in H^2(\TT)\\ \|\vp\|_{H^2}\le 1
\end{array}
}\bigg{|}\langle\pat u_n(t),\vp\rangle\bigg{|}.
\]
Then, we estimate as
\begin{align*}
\av{\int_{\TT^2}\pat u^{(n)}\vp\,dx}&\le c\int_{\TT^2}\pare{\av{u^{(n)}}+\av{\D u^{(n)}}+(\D u^{(n)})^2}\av{\D \vp}\,dx\\
&\q+c\int_{\TT^2}\av{\det D^2u^{(n)}}\vp\,dx\\
&\le c\pare{1+\norm{u^{(n)}}_{W^{2,\infty}}}\norm{u^{(n)}}_{H^{2}}\norm{\vp}_{H^2}.
\end{align*}
Then,
\[
\norm{\pat u^{(n)}(t)}_{H^{-2}}\le c\pare{1+\norm{u^{(n)}}_{W^{2,\infty}}}\norm{u^{(n)}}_{H^{2}},
\]
 and
\[
\norm{\pat u^{(n)}(t)}_{L^2(H^{-2})}\le c\pare{1+\norm{u^{(n)}}_{L^{\infty}(W^{2,\infty})}}\norm{u^{(n)}}_{L^2(H^{2})}<c.
\]

\textbf{Passing to the limit:} We want to take the limit in $n$ in
\begin{align*}
&\int_{\TT^2} P_n u_{0}\vp(0)\,dx+\iint_{\TT^2\times(0,T)} u^{(n)}\pat\vp+\pare{K_0 u^{(n)}+\frac{K_3}{2}(\D u^{(n)})^2}\D\vp\,dx\,dt\\
& \qq +\iint_{\TT^2\times(0,T)}  2K_1 \vp \det D^2u^{(n)}
-K_2u^{(n)}\D^2\vp  \,dx\,dt=0,
\end{align*}
 being $\vp\in W^{1,1}(0,T;L^1(\TT^2))\cap L^1(0,T;W^{4,1}(\TT^2))\cap L^2(0,T;H^2(\TT^2))$.

We only detail the convergence of
\[
\iint_{\TT^2\times(0,T)}(\D u^{(n)})^2\D\vp  \,dx\,dt
\]
because the one of
\[
\iint_{\TT^2\times(0,T)}\vp \det D^2u^{(n)} \,dx\,dt
\]
is similar to the first one, and the others follow from the assumptions on $\vp$ and the weak-$*$ convergence \eqref{eq:weak*1-hessian}. We have
\begin{align*}
&\iint_{\TT^2\times(0,T)}\pare{(\D u^{(n)})^2-(\D u)^2}\D\vp  \,dx\,dt\\
&=\iint_{\TT^2\times(0,T)}\D(u^{(n)}-u)\D(u^{(n)}+u)\D\vp  \,dx\,dt\\
&\le \norm{u^{(n)}+u}_{L^\infty(W^{2,\infty})}\iint_{\TT^2\times(0,T)}\av{\D(u^{(n)}-u)}\av{\D\vp}  \,dx\,dt
\end{align*}
thanks to \eqref{eq:weak*1-hessian}. We now apply H\"older`s inequality, obtaining that
\begin{align*}
\iint_{\TT^2\times(0,T)}\pare{(\D u_n)^2-(\D u)^2}\D\vp  \,dx\,dt
&\le \norm{u_n+u}_{L^\infty(W^{2,\infty})}\norm{u_n-u}_{L^2(H^{2})}\norm{\vp}_{L^2(H^{2})},
\end{align*}
which converges to zero as $n\to\infty$ by \eqref{eq:strong-hessian}.

\section{Proof of Theorem \ref{teo:ex-porous}}
\textbf{Approximate problem:} We Observe that in the new variable $v$, problem \eqref{eq:pb-porous-thin-v} is equivalent to
\begin{equation*}
\begin{cases}
\begin{array}{ll}
\ds\pat v
=-v,_{iijj}
-v,_iv,_{jji}-v v,_{iijj}-\chi p(p-1) (1+v)^{p-2} v,_iv,_i -\chi p (1+v)^{p-1}v,_{ii}
 & \text{in}\q  (0,T)\times\TT^2, \\
\ds v^{(n)}(0,x)=u_0(x)-1 &\t{in}\q \TT^2.
\end{array}
\end{cases}
\end{equation*}
We define the following approximate problem
\begin{equation}\label{eq:pb-porous-thin-v-einstein-appr}
\begin{cases}
\begin{array}{lll}
&\ds\pat v^{(n)}
=-v^{(n)},_{iijj}
-v^{(n)},_iv^{(n)},_{jji}-v^{(n)} v^{(n)},_{iijj}&\\
&\qq\qq-\chi p(p-1) (1+v^{(n)})^{p-2} v^{(n)},_iv^{(n)},_i -\chi p (1+v^{(n)})^{p-1}v^{(n)},_{ii}
 & \text{in}\q  (0,T)\times\TT^2, \\
&\ds v(0,x)=P_n (u_0(x)-1) &\t{in}\q \TT^2.
\end{array}
\end{cases}
\end{equation}

\textbf{A priori estimates in Wiener spaces:} Let us omit the superscript $(n)$ in the following computations. We use that
\[
(a+b)^n=\sum_{k=0}^{n}\binom{n}{k}a^{n-k}b^k
\]
and we rewrite each term of the r.h.s. of \eqref{eq:pb-porous-thin-v-einstein-appr} in Fourier:
\begin{align*}
\w{v,_{iijj}}(t,k)&=|k|^4\w{u}(t,k),\\
\w{v,_iv,_{jji}}(t,k)&=\sum_{m\in\ZZ^2}m_i(k_i-m_i)|k-m|^2\w{v}(t,m)\w{v}(t,k-m),\\
\w{vv,_{iijj}}(t,k)&=\sum_{m\in\ZZ^2}|k-m|^4\w{v}(t,m)\w{v}(t,k-m),
\end{align*}
\begin{align*}
&\w{(1+v)^{p-2}v,_iv,_i}(t,k)\\
&=\sum_{q=0}^{p-2}\binom{p-2}{q}\w{v^qv,_iv,_i}(t,k)\\
%
&=-\sum_{q=0}^{p-2}\binom{p-2}{q}\sum_{m^1\in\ZZ^2}\ldots\sum_{m^{q+1}\in\ZZ^2}(k_i-m_i^1)(m_i^1-m_i^2)\w{v}(t,k-m^1)\w{v}(t,m^1-m^2)\times\\
&\qq\times\prod_{\ell=2}^{q}\w{v}(t,m^\ell-m^{\ell+1})\w{v}(t,m^{q+1})\\
&=-\sum_{q=0}^{p-2}\frac{(p-2)!}{q!(p-2-q)!}\m{N}_1^{(q)}(t,k), \\
&\w{(1+v)^{p-1}v,_{ii}}(k,t)\\
&=\sum_{q=0}^{p-1}\binom{p-1}{q}\w{v^qv,_{ii}}(t,k)\\
%
%
&
=-\sum_{q=0}^{p-1}\binom{p-1}{q}\sum_{m^1\in\ZZ^2}\ldots\sum_{m^{q}\in\ZZ^2}\av{k-m^1}^2\w{v}(t,k-m^1) \prod_{\ell=1}^{q-1}\w{v}(t,m^\ell-m^{\ell+1})\w{v}(t,m^{q})\\
&=-\sum_{q=0}^{p-1}\frac{(p-1)!}{q!(p-1-q)!}\m{N}_2^{(q)}(t,k).
\end{align*}
Then, the Fourier coefficients of \eqref{eq:pb-porous-thin-v-einstein-appr} satisfy
\begin{align*}
\pat\w{v}(t,k)&=(-|k|^4+\chi|k|^2)\w{v}(t,k)-\sum_{m\in\ZZ^2}m_i(k_i-m_i)|k-m|^2\w{v}(t,m)\w{v}(t,k-m)\\
&\q-\sum_{m\in\ZZ^2}|k-m|^4\w{v}(t,m)\w{v}(t,k-m)
%
%
 +\chi p !\sum_{q=0}^{p-2}\frac{\m{N}_1^{(q)}(t,k)}{q!(p-2-q)!}  +\chi p! \sum_{q=1}^{p-1}\frac{\m{N}_2^{(q)}(t,k)}{q!(p-1-q)!}.
\end{align*}
We recall \eqref{eq:der-t} to deduce
\[
\sum_{k\in\ZZ^2}	\pat\av{ \w{v}(k,t)}=\frac{d}{dt}\norm{v(t)}_{A^0},\q
\sum_{k\in\ZZ^2}|k|^4\av{\w{v}(k,t)}=\norm{v(t)}_{A^4},
\]
\begin{align*}
\sum_{k\in\ZZ^2}\av{\sum_{m\in\ZZ^2}m_i(k_i-m_i)|k-m|^2\w{v}(t,m)\w{v}(t,k-m)}&\le  \norm{v(t)}_{A^1}\norm{v(t)}_{A^3}
\\
&\le  \norm{v(t)}_{A^0}\norm{v(t)}_{A^4},\\
%
\sum_{k\in\ZZ^2}\av{\sum_{m\in\ZZ^2}|k-m|^4\w{v}(t,m)\w{v}(t,k-m)}&\le  \norm{v(t)}_{A^0}\norm{v(t)}_{A^4},
\end{align*}
\begin{align*}
\sum_{q=0}^{p-2}\frac{1}{q!(p-2-q)!}\sum_{k\in\ZZ^2}\av{\m{N}_1^{(q)}(t,k)}
&\le c \sum_{q=0}^{p-2}\norm{v(t)}_{A^0}^q\norm{v(t)}_{A^1}^2
\le c \sum_{q=0}^{p-2}\norm{v(t)}_{A^0}^{q+1}\norm{v(t)}_{A^2} ,\\
\sum_{q=0}^{p-1}\frac{1}{q!(p-1-q)!}\sum_{k\in\ZZ^2}\av{\m{N}_2^{(q)}(t,k)}
&\le c\sum_{q=0}^{p-1}\norm{v(t)}_{A^0}^q\norm{v(t)}_{A^2},
%
%
\end{align*}
and Poincar\'e inequality for Wiener spaces, we estimate the $A^0$ seminorm of $v$ as
\begin{align*}
\frac{d}{dt}\norm{v(t)}_{A^0}&\le -\pare{(1-\chi)-2\norm{v(t)}_{A^0}}\norm{v(t)}_{A^4}+ c \chi p!\pare{
\sum_{q=0}^{p-2}\norm{v(t)}_{A^0}^{q+1}+
 \sum_{q=0}^{p-1}\norm{v(t)}_{A^0}^q
}\norm{v(t)}_{A^2}
\\
&=-\pare{1-\chi-2\norm{v(t)}_{A^0}}\norm{v(t)}_{A^4} +c\chi p!\pare{\norm{v(t)}_{A^0}+2
\sum_{q=1}^{p-1}\norm{v(t)}_{A^0}^q
}\norm{v(t)}_{A^2}.
\end{align*}
We estimate the $\norm{v(t)}_{A^2}$ as
\[
\norm{v(t)}_{A^2}\le \norm{v(t)}_{A^0}^\frac{1}{2}\norm{v(t)}_{A^4}^\frac{1}{2}\le \frac{1}{2}\norm{v(t)}_{A^0}  + \frac{1}{2}\norm{v(t)}_{A^4},
\]
obtaining that
\begin{align*}
\frac{d}{dt}\norm{v(t)}_{A^0}&\le -
\pare{
1-\chi-2\norm{v(t)}_{A^0}
-\frac{c\chi p!}{2}\pare{\norm{v(t)}_{A^0}+2
\sum_{q=1}^{p-1}\norm{v(t)}_{A^0}^q
}
}\norm{v(t)}_{A^4} \\
&\q
+\frac{c\chi p!}{2}\pare{\norm{v(t)}_{A^0}+2
\sum_{q=1}^{p-1}\norm{v(t)}_{A^0}^q
}\norm{v(t)}_{A^0}.
\end{align*}
The smallness condition
\begin{equation*}
1-\chi-2\norm{v_0}_{A^0}
-\frac{c\chi p!}{2}\pare{\norm{v_0}_{A^0}+2
\sum_{q=1}^{p-1}\norm{v_0}_{A^0}^q
}>0,
\end{equation*}
implies that, for small times,
\begin{align*}
&\frac{d}{dt}\norm{v(t)}_{A^0}+
\pare{
1-\chi-2\norm{v_0}_{A^0}
-\frac{c\chi p!}{2}\pare{\norm{v_0}_{A^0}+2
\sum_{q=1}^{p-1}\norm{v_0}_{A^0}^q
}
}\norm{v(t)}_{A^4} \\
&\le
\frac{c\chi p!}{2}\pare{\norm{v(t)}_{A^0}+2
\sum_{q=1}^{p-1}\norm{v(t)}_{A^0}^q
}\norm{v(t)}_{A^0}\\
&\le
\frac{c\chi p!}{2}\pare{\norm{v_0}_{A^0}+2
\sum_{q=1}^{p-1}\norm{v_0}_{A^0}^q
}\norm{v_0}_{A^0}.
\end{align*}
We extend this inequality to all times using a  contradiction argument in time.\\
As a consequence, we find the uniform boundedness in the following space
\begin{equation*}
 W^{1,1}(0,T;A^0(\TT^2))\cap L^1(0,T;A^4(\TT^2)).
\end{equation*}
Invoking a Poincar\'e inequality in Wiener spaces we conclude the decay estimate
\begin{equation*}
\norm{v^{(n)}(t)}_{A^0}\le \exp
\pare{-\pare{
1-\chi-2\norm{v_0}_{A^0}
-c\chi p!\pare{\norm{v_0}_{A^0}+2
\sum_{q=1}^{p-1}\norm{v_0}_{A^0}^q
}}
t}
\norm{v_0}_{A^0},
\end{equation*}
for $\chi$ eventually smaller.

\textbf{Compactness results:} Reasoning as for the previous problem \eqref{eq:pb-hessian-einstein-appr},  we have that
\begin{align*}
v^{(n)}&\to v && \t{a.e. }(0,T)\times \TT^2,\label{eq:ae-porous-thin}\\
v^{(n)}&\overset{*}{\rightharpoonup} v && \t{in } L^{\infty}((0,T)\times\TT^2),
\\
v^{(n)}&\overset{*}{\rightharpoonup} v && \t{in } \m{M}(0,T;W^{4,\infty}( \TT^2)),
\\
v^{(n)}&\overset{*}{\rightharpoonup} v && \t{in } L^{\frac{4}{m}}(0,T;W^{m,\infty}( \TT^2))\q\t{with}\q 0< m< 4,
\\
%
%
v^{(n)}&\rightharpoonup v && \t{in } L^{2}(0,T;H^{2}( \TT^2)),
\\
%
%
v^{(n)}&\to v && \t{in } L^{2}(0,T;H^{r}( \TT^2))\q\t{with}\q 0\le r<2,
\\
%
%
v^{(n)}&\to v && \t{in } L^{q}(0,T;H^{2}( \TT^2))\q\t{with}\q 1\le q<2
\end{align*}
up to subsequences.

\textbf{Passing to the limit:} We want to take the limit in $n$ in
\begin{align*}
&-\int_{\TT^2} P_n v_0\vp(0)\,dx\\
&\q +\iint_{\TT^2\times(0,T)} -v^{(n)}\pat\vp-v^{(n)}\nabla\D v^{(n)}\cdot \nabla \vp+\chi (v^{(n)}+1)^{p-1}\D\vp  \,dx\,dt=0,
\end{align*}
for every $\vp\in W^{1,1}(0,T;L^1(\TT^2))\cap L^{q}(0,T;H^2(\TT^2))$ and $1\leq q<2$. Then, using that, due to interpolation we have strong convergence in
$$
L^r(0,T;H^3),
$$
and since
\[
(v_n+1)^{p-1}-(v+1)^{p-1}=\sum_{q=1}^{p-1}v_n^q-v^q=(v_n-v)\sum_{q=2}^{p-1}v_n^qv^{p-1-q}
\]
we deduce that
\begin{align*}
&\av{\iint_{\TT^2\times(0,T)}\pare{(v_n+1)^{p-1}-(v+1)^{p-1}}\D \vp\,dx\,dt}\\&\le \sum_{q=2}^{p-1}\norm{v_n}_{L^\infty(L^\infty)}^q\norm{v}_{L^\infty(L^\infty)}^{p-1-q}\norm{v_n-v}_{L^2(L^2)}\norm{\vp}_{L^{2}(H^2)}\to 0,
\end{align*}
and we can pass to the limit.

\section{Conclusion}
In this paper we have studied two fourth order nonlinear parabolic equations. These equations arise in the study of epitaxial growth and in thin films. Such quasilinear PDEs are mathematically challenging due to the high number of derivatives in the main part. For these PDEs we have established a number of global in time existence results in the non-standard Wiener spaces. These functional spaces allow us to take full advantage of the parabolic structure.

In particular, one of our main contributions has been the improvement of the previous global in time result contained in \cite{EGP}. Indeed, the authors in \cite{EGP} establish a global in time results for initial data with small energy akin to the $H^2$ norm. However, with our techniques we can prove a global in time result imposing a size condition in the Wiener algebra $A^0$. The Wiener algebra shares the same scaling as $L^\infty$ and both are critical spaces with respect to the scaling of the equation.

\section*{Acknowledgments}

R.G-B was supported by the project "Mathematical Analysis of Fluids and Applications" Grant PID2019-109348GA-I00 funded by MCIN/AEI/ 10.13039/501100011033 and acronym "MAFyA". This publication is part of the project PID2019-109348GA-I00 funded by MCIN/ AEI /10.13039/501100011033. R.G-B is also supported by a 2021 Leonardo Grant for Researchers and Cultural Creators, BBVA Foundation. The BBVA Foundation accepts no responsibility for the opinions, statements, and contents included in the project and/or the results thereof, which are entirely the responsibility of the authors.\\
 The work of M.M. was partially supported by Grant RYC2021-033698-I, funded by the Ministry of Science and Innovation/State Research Agency/10.13039/501100011033 and by the European Union "NextGenerationEU/Recovery, Transformation and Resilience Plan".\\
Both authors are funded by  the project "An\'alisis Matem\'atico Aplicado y Ecuaciones Diferenciales" Grant PID2022-141187NB-I00 funded by MCIN /AEI /10.13039/501100011033 / FEDER, UE and acronym "AMAED". This publication is part of the project PID2022-141187NB-I00 funded by MCIN/ AEI /10.13039/501100011033.

\end{document}